\documentclass[11pt,twoside]{article}
\usepackage{amsmath, amsthm, amscd, amsfonts, amssymb, graphicx, color}

\date{}

\begin{document}

\centerline{}

\centerline {\Large{\bf  Generalized $p$-fusion frame in separable Banach space}}

\newcommand{\mvec}[1]{\mbox{\bfseries\itshape #1}}
\centerline{}
\centerline{\textbf{Prasenjit Ghosh}}
\centerline{Department of Pure Mathematics, University of Calcutta,}
\centerline{35, Ballygunge Circular Road, Kolkata, 700019, West Bengal, India}
\centerline{e-mail: prasenjitpuremath@gmail.com}
\centerline{}
\centerline{\textbf{T. K. Samanta}}
\centerline{Department of Mathematics, Uluberia College,}
\centerline{Uluberia, Howrah, 711315,  West Bengal, India}
\centerline{e-mail: mumpu$_{-}$tapas5@yahoo.co.in}

\newtheorem{Theorem}{\quad Theorem}[section]

\newtheorem{definition}[Theorem]{\quad Definition}

\newtheorem{theorem}[Theorem]{\quad Theorem}

\newtheorem{remark}[Theorem]{\quad Remark}

\newtheorem{corollary}[Theorem]{\quad Corollary}

\newtheorem{note}[Theorem]{\quad Note}

\newtheorem{lemma}[Theorem]{\quad Lemma}

\newtheorem{example}[Theorem]{\quad Example}

\newtheorem{result}[Theorem]{\quad Result}
\newtheorem{conclusion}[Theorem]{\quad Conclusion}

\newtheorem{proposition}[Theorem]{\quad Proposition}

\begin{abstract}
\textbf{\emph{Concepts of $g$-fusion frame and $gf$-Riesz basis in a Hilbert to a Banach space is being presented.\,Some properties of $g$-fusion frame and $gf$-Riesz basis in Banach space have been developed.\,We discuss perturbation results of $g$-fusion frame in a Banach space.\,Finally, we construct $g$-$p$-fusion frames in Cartesian product of Banach spaces and tensor product of Banach spaces.}}
\end{abstract}


\section{Introduction}

In recent times, several generalization of frame for separable Hilbert space have been introduced.\,Some of them are \,$K$-frame \cite{Gavruta,Xiao}, \,$g$-frame \cite{Sun}, fusion frame \cite{Kutyniok} and so on.\,The combination of \,$g$-frame and fusion frame is known as generalized fusion frame or \,$g$-fusion frame.\,Sadri et al.\,\cite{Ahmadi} presented \,$g$-fusion frame to generalize the theory of fusion frame and \,$g$-frame.\,P.\,Ghosh and T.\,K.\,Samanta  \cite{P} presented the stability of dual \,$g$-fusion frames in Hilbert spaces.\,These frames were further studied in \cite{Ghosh, GG}.

A. Aldroubi et al.\,\cite{A} introduced \,$p$-frame in a Banach space and discussed some of its properties.\,Chistensen and stoeva \cite{O} also developed \,$p$-frame in separable Banach space.\,M. R. Abdollahpour et al. \cite{M} introduced the \,$p\,g$-frames in Banach spaces.\,The generalization of the \,$g$-frame and \,$g$-Riesz Basis in a complex Hilbert space to a complex Banach space was also studied by Xiang-Chun Xio et al.\,\cite{X}.

In this paper, we generalize the notion of \,$g$-fusion frame in a Hilbert space to a Banach space and establish some of its properties.\,Generalized Riesz basis in Banach space is also discussed.\;The relation between \,$g$-$p$-fusion frame and \,$q$-$g\,f$-Riesz basis is obtained.\,We describe some perturbation results of \,$g$-$p$-fusion frame in Banach space.\,At the end, we present \,$g$-$p$-fusion frame in tensor product of Banach spaces.

\section{Preliminaries}

\smallskip\hspace{.6 cm} Throughout this paper,\;$X$\, is considered to be a separable Banach space over the field \,$\mathbb{K}\,(\,\mathbb{R}\, \,\text{or}\; \,\mathbb{C}\,)$\, and \,$X^{\,\ast}$, its dual space.\,$I,\, J$\, denotes the subset of  natural numbers \,$\mathbb{N}$.\,$\left\{\,X_{i}\,\right\}_{i \,\in\, I}$\, is a sequence of Banach spaces and \,$\left\{\,V_{i}\,\right\}_{ i \,\in\, I}$\, is a collection of closed subspaces of \,$X$.\,$\mathcal{B}\,(\,X,\, X_{i}\,)$\, are the collection of all bounded linear operators from \,$X \;\text{to}\; X_{i}$\, and in particular, \,$\mathcal{B}\,(\,X\,)$\, denotes the space of all bounded linear operators on \,$X$.\,It is assumed that \,$p \,\in\, (\,1,\,\infty\,)$\, and when \,$p$\, and \,$q$\, are used in a same assertion, they satisfy the relation \,$1 \,/\, p \,+\, 1 \,/\, q \,=\, 1$.

\begin{theorem}\label{thm1}\cite{H}
If \,$U \,:\, X \,\to\, Y$\, is a bounded operator from a Banach space \,$X$\, into a Banach space \,$Y$\, then its adjoint \,$U^{\,\ast} \,:\, Y^{\,\ast} \,\to\, X^{\,\ast}$\, is surjective if and only if \,$U$\, has a bounded inverse on \,$\mathcal{R}_{U}$(\,range of $U$\,).
\end{theorem}

\begin{definition}
Let \,$1 \,<\, p \,<\, \infty$.\,A countable family \,$\left\{\,g_{\,i}\,\right\}_{i \,\in\, I} \,\subset\, X^{\,\ast}$\, is said to be a p-frame for \,$X$\, if there exist constants \,$0 \,<\, A \,\leq\, B \,<\, \infty$\, such that
\[A \,\left\|\,f\,\right\|_{X} \,\leq\, \left(\,\sum\limits_{\,i \,\in\, I}\,\left|\,g_{\,i}\,(\,f\,)\,\right|^{\,p}\,\right)^{1 \,/\, p} \,\leq\, B \,\left\|\, f\,\right\|_{X}\; \;\forall\; f \,\in\, X.\]
\end{definition}

\begin{definition}\cite{X}
A sequence \,$\left\{\,\Lambda_{i} \,\in\, \mathcal{B}\,(\,X,\, X_{i}\,) \,:\, i \,\in\, I\,\right\}$\, is called a generalized   p-frame or g-p-frame for \,$X$\; with respect to \,$\left\{\,X_{i}\,\right\}_{i \,\in\, I}$\; if there exist two positive constants \,$A$\; and \,$B$\; such that
\[A \;\left \|\, f \,\right \|_{X} \,\leq\, \left(\,\sum\limits_{\,i \,\in\, I}\, \left\|\,\Lambda_{i}\,(\,f\,) \,\right\|^{\,p}\,\right)^{1 \,/\, p} \,\leq\, B \; \left\|\, f \, \right\|_{X}\; \;\forall\; f \,\in\, X.\]
\,$A$\; and \,$B$\; are called the lower and upper frame bounds, respectively.
\end{definition}

\begin{definition}\cite{X}
Define the linear space
\[l^{\,p}\left(\,\left\{\,X_{i}\,\right\}_{i \,\in\, I}\,\right) \,=\, \left\{\,\{\,f_{\,i}\,\}_{i \,\in\, I} \,:\, f_{\,i} \,\in\, X_{i},\; \sum\limits_{\,i \,\in\, I}\, \left\|\,f_{\,i}\,\right \|^{\,p} \,<\, \infty \,\right\}.\]
Then it is a complex Banach space with respect to the norm is defined by
\[\left\|\,\left\{\,f_{\,i}\,\right\}_{i \,\in\, I}\,\right\| \,=\, \left(\,\sum\limits_{\,i \,\in\, I}\, \left\|\,f_{\,i}\,\right \|^{\,p}\,\right)^{1 \,/\, p}\]
\end{definition}

\begin{lemma}\cite{X}
Let \,$p \,>\, 1,\, q \,>\, 1$\, be such that \,$1 \,/\, p \,+\, 1 \,/\, q \,=\, 1$.\,Then the adjoint space of \,$l^{\,p}\left(\,\left\{\,X_{i}\,\right\}_{i \,\in\, I}\,\right)$\, is \,$l^{\,q}\left(\,\left\{\,X^{\,\ast}_{i}\,\right\}_{i \,\in\, I}\,\right)$, where \,$X^{\,\ast}_{i}$\, is the adjoint space of \,$X_{i}$\, for \,$i \,\in\, I$.
\end{lemma}

\begin{definition}\cite{M,X}
Let \,$\left\{\,\Lambda_{i} \,\in\, \mathcal{B}\,(\,X,\, X_{i}\,) \,:\, i \,\in\, I\,\right\}$\, be a generalized p-frame or g-p-frame for \,$X$.\,Then the operator defined by
\[U \,:\, X \,\to\, l^{\,p}\,\left(\,\left\{\,X_{i}\,\right\}_{i \,\in\, I}\,\right),\; \;U\,f \,=\, \left\{\,\Lambda_{i}\,(\,f\,)\,\right\}_{i \,\in\, I}\; \;\forall\; f \,\in\, X.\]
is called the analysis operator and the operator given by 
\[T \,:\, l^{\,q}\,\left(\,\left\{\,X^{\,\ast}_{i}\,\right\}_{i \,\in\, I}\,\right) \,\to\, X^{\,\ast}\] 
\[T\,\left(\,\left\{\,g_{i}\,\right\}_{i \,\in\, I}\,\right) \,=\, \sum\limits_{i \,\in\, I}\,\Lambda^{\,\ast}_{i}\,g_{\,i}\; \;\forall\, \left\{\,g_{i}\,\right\}_{i \,\in\, I} \,\in\, l^{\,q}\,\left(\,\left\{\,X^{\,\ast}_{i}\,\right\}_{i \,\in\, I}\,\right)\]
is called synthesis operator.  
\end{definition}

\begin{definition}\cite{Ahmadi}
Let \,$\left\{\,v_{i}\,\right\}_{ i \,\in\, I}$\, be a collection of positive weights and \,$\left\{\,H_{i}\,\right\}_{i \,\in\, I}$\, be a collections of Hilbert spaces and \,$\left\{\,V_{i}\,\right\}_{i \,\in\, I}$\, be a family of closed subspaces of a Hilbert space \,$H$.\;Then the family \,$\Lambda \,=\, \{\,\left(\,V_{i},\, \Lambda_{i},\, v_{i}\,\right)\,\}_{i \,\in\, I}$\; is called a generalized fusion frame or a g-fusion frame for \,$H$\; respect to \,$\left\{\,H_{i}\,\right\}_{i \,\in\, I}$\; if there exist constants \,$0 \,<\, A \,\leq\, B \,<\, \infty$\; such that
\begin{equation}\label{eqq1}
A \;\left \|\,f \,\right \|^{\,2} \,\leq\, \sum\limits_{\,i \,\in\, I}\,v_{i}^{\,2}\, \left\|\,\Lambda_{i}\,P^{\,\prime}_{\,V_{i}}\,(\,f\,) \,\right\|^{\,2} \,\leq\, B \; \left\|\, f \, \right\|^{\,2}\; \;\forall\; f \,\in\, H,
\end{equation}
where \,$P^{\,\prime}_{\,V_{i}}$\, is the orthogonal projection of \,$H$\, onto \,$V_{i}$.\,The constants \,$A$\; and \,$B$\; are called the lower and upper bounds of g-fusion frame, respectively.\,If \,$\Lambda$\; satisfies the right inequality of (\ref{eqq1}), it is called a g-fusion Bessel sequence with bound \,$B$\; in \,$H$. 
\end{definition}

\section{$g$-$p$-fusion frame and it's properties}

\smallskip\hspace{.6 cm} In this section, we develop the generalized fusion frame and generalized Riesz basis for Banach space.    

\begin{definition}\label{deff1}
Let \,$p \,>\, 1$\, and \,$\left\{\,v_{i}\,\right\}_{ i \,\in\, I}$\, be a collection of positive weights i\,.\,e., \,$v_{\,i} \,>\, 0$.\,Let \,$\Lambda_{i} \,\in\,  \mathcal{B}\,(\,X,\, X_{i}\,)$\, and \,$\left\{\,P_{\,V_{i}}\,\right\}$\, be non-trivial linear projections of \,$X$\, onto \,$V_{i}$\, such that \,$P_{\,V_{i}}\,(\,X\,) \,=\, V_{i}$, for each \,$i \,\in\, I$.\,Then the family \,$\Lambda \,=\, \left\{\,\left(\,V_{i},\, \Lambda_{i},\, v_{i}\,\right)\,\right\}_{i \,\in\, I}$\, is called a generalized \,$p$-fusion frame or a g-p-fusion frame for \,$X$\, with respect to \,$\left\{\,X_{i}\,\right\}_{i \,\in\, I}$\, if there exist constants \,$0 \,<\, A \,\leq\, B \,<\, \infty$\; such that
\begin{equation}\label{eq1}
A \,\left \|\,f \,\right \| \,\leq\, \left(\,\sum\limits_{\,i \,\in\, I}\,v_{i}^{\,p}\,\left\|\,\Lambda_{i}\,P_{\,V_{i}}\,(\,f\,) \,\right\|^{\,p}\,\right)^{1 \,/\, p} \,\leq\, B \,\left\|\,f\, \right\|\, \,\forall\, f \,\in\, X.
\end{equation}
The constants \,$A$\, and \,$B$\, are called the lower and upper bounds of g-p-fusion frame, respectively.\,If \,$A \,=\, B$\; then \,$\Lambda$\; is called tight g-p-fusion frame and if \;$A \,=\, B \,=\, 1$\, then we say \,$\Lambda$\; is a Parseval g-p-fusion frame.\;If  \,$\Lambda$\; satisfies only the right inequality of (\ref{eq1}), it is called a g-p-fusion Bessel sequence with bound \,$B$\; in \,$X$.   
\end{definition}

\begin{remark}
Suppose that \,$\Lambda \,=\, \left\{\,\left(\,V_{i},\, \Lambda_{i},\, v_{i}\,\right)\,\right\}_{i \,\in\, I}$\, is a tight g-p-fusion frame for \,$X$\, with bound \,$A$.\,Then for all \,$f \,\in\, X$, we have
\[\left(\,\sum\limits_{\,i \,\in\, I}\,v_{i}^{\,p}\, \left\|\,\Lambda_{i}\,P_{\,V_{i}}\,(\,f\,) \,\right\|^{\,p}\,\right)^{1 \,/\, p} = A\,\left\|\,f\, \right\|\]
\[ \,\Rightarrow\, \left(\,\sum\limits_{i \,\in\, I}\,v_{i}^{\,p}\, \left\|\,A^{\,-\, 1 }\,\Lambda_{j}\,P_{\,V_{i}}\,(\,f\,) \,\right\|^{\,p}\,\right)^{1 \,/\, p} = \left\|\,f\, \right\|.\]
This verify that \,$\left\{\,\left(\,V_{i},\, A^{\,-\, 1 }\,\Lambda_{i},\, v_{i}\,\right)\,\right\}_{i \,\in\, I}$\, is a Parseval g-p-fusion frame for \,$X$.   
\end{remark}

\begin{theorem}\label{th1.20}
Let \,$\Lambda$\, be a g-p-fusion frame for \,$X$\, with respect to \,$\left\{\,X_{i}\,\right\}_{i \,\in\, I}$\, having bounds \,$A,\,B$.\,Suppose \,$U \,\in\, \mathcal{B}\,(\,X\,)$\, be an invertible operator on \,$X$.\,Then \,$\Gamma \,=\, \left\{\,\left(\,U\,V_{i},\, \Lambda_{i}\,P_{\,V_{i}}\,U,\, v_{i}\,\right)\,\right\}_{i \,\in\, I}$\, is a g-p-fusion frame for \,$X$, provided \,$P_{\,V_{i}}\,U\,P_{\,U\,V_{i}} \,=\, P_{\,V_{i}}\,U$, for \,$i \,\in\, I$. 
\end{theorem}

\begin{proof}
For each \,$f \,\in\, X$, we have
\begin{align*}
&\left(\,\sum\limits_{\,i \,\in\, I}\,v_{i}^{\,p}\, \left\|\,\Lambda_{i}\,P_{\,V_{i}}\,U\,P_{\,U\,V_{i}}\,(\,f\,) \,\right\|^{\,p}\,\right)^{1 \,/\, p} \,=\, \left(\,\sum\limits_{\,i \,\in\, I}\,v_{i}^{\,p}\, \left\|\,\Lambda_{i}\,P_{\,V_{i}}\,U\,(\,f\,)\,\right\|^{\,p}\,\right)^{1 \,/\, p}\\
&\leq\, B\,\left\|\,U\,f\,\right\| \,\leq\, B\,\left\|\,U\,\right\|\,\left\|\,f\, \right\|\; \;[\;\text{since $\Lambda$\, is a $g$-$p$-fusion frame}\;].
\end{align*}
On the other hand
\begin{align*}
&\left(\,\sum\limits_{\,i \,\in\, I}\,v_{i}^{\,p}\, \left\|\,\Lambda_{j}\,P_{\,V_{i}}\,U\,P_{\,U\,V_{i}}\,(\,f\,) \,\right\|^{\,p}\,\right)^{1 \,/\, p} \,=\, \left(\,\sum\limits_{\,i \,\in\, I}\,v_{i}^{\,p}\, \left\|\,\Lambda_{i}\,P_{\,V_{i}}\,U\,(\,f\,)\,\right\|^{\,p}\,\right)^{1 \,/\, p}\\
&\geq\, A\,\left\|\,U\,f\,\right\| \,\geq\, A\,\left\|\,U^{\,-\, 1}\,\right\|^{\,-\, 1}\,\left\|\,f\, \right\|\; \;[\;\text{since $U$ is invertible}\;].
\end{align*}
Hence, \,$\Gamma$\, is a \,$g$-$p$-fusion frame for \,$X$\, with bounds \,$B\,\left\|\,U\,\right\|$\, and \,$A\,\left\|\,U^{\,-\, 1}\,\right\|^{\,-\, 1}$.
\end{proof}

\begin{theorem}
Let \,$\Lambda$\, be a g-p-fusion frame for \,$X$\, with respect to \,$\left\{\,X_{i}\,\right\}_{i \,\in\, I}$\, having bounds \,$A,\,B$\, and \,$U \,:\, X \,\to\, X$\, be a bounded linear operator such that for each \,$i \,\in\, I,\, \,P_{\,V_{i}}\,U\,P_{\,U\,V_{i}} \,=\, P_{\,V_{i}}\,U$.\,Then the family \,$\Gamma \,=\, \left\{\,\left(\,U\,V_{i},\, \Lambda_{i}\,P_{\,V_{i}}\,U,\, v_{i}\,\right)\,\right\}_{i \,\in\, I}$\, is a g-p-fusion frame for \,$X$\, if and only if \,$U$\, is bounded below. 
\end{theorem}

\begin{proof}
Let \,$\Gamma$\, be a \,$g$-$p$-fusion frame for \,$X$\, with bounds \,$C$\, and \,$D$.\,Then
\[C\,\left\|\,f\,\right\| \,\leq\, \left(\,\sum\limits_{\,i \,\in\, I}\,v_{i}^{\,p}\, \left\|\,\Lambda_{i}\,P_{\,V_{i}}\,U\,P_{\,U\,V_{i}}\,(\,f\,) \,\right\|^{\,p}\,\right)^{1 \,/\, p} \,\leq\, D \,\left\|\,f\, \right\|\, \,\forall\ f \,\in\, X.\]
\begin{equation}\label{eq1.1}
\Rightarrow\, C \,\left \|\,f \,\right \| \,\leq\, \left(\,\sum\limits_{\,i \,\in\, I}\,v_{i}^{\,p}\, \left\|\,\Lambda_{i}\,P_{\,V_{i}}\,U\,(\,f\,) \,\right\|^{\,p}\,\right)^{1 \,/\, p} \,\leq\, D \,\left\|\,f\, \right\|
\end{equation}
Since \,$\Lambda$\, is a \,$g$-$p$-fusion frame with bounds \,$A$\, and \,$B$, in (\ref{eq1}), replacing \,$f$\, by \,$U\,f$, we get
\begin{equation}\label{eq1.2}
A \,\left \|\,U\,f\,\right \| \,\leq\, \left(\,\sum\limits_{\,i \,\in\, I}\,v_{i}^{\,p}\, \left\|\,\Lambda_{i}\,P_{\,V_{i}}\,U\,(\,f\,)\,\right\|^{\,p}\,\right)^{1 \,/\, p} \,\leq\, B \,\left\|\,U\,f\,\right\|\, \,\forall\ f \,\in\, X.
\end{equation}
Now, from (\ref{eq1.1}) and (\ref{eq1.2}), for each \,$f \,\in\, X$, we can write
\[C\,\left\|\,f\,\right\| \,\leq\, B \,\left\|\,U\,f\,\right\| \,\Rightarrow\, \left\|\,U\,f\,\right\|  \,\geq\, \dfrac{C}{B}\,\left\|\,f\,\right\|.\]
This shows that \,$U$\, is bounded below.\\

Conversely, suppose that there exists \,$M \,>\, 0$\, such that \,$\left\|\,U\,f\,\right\|  \,\geq\, M\,\left\|\,f\,\right\|$.\,Now, for each \,$f \,\in\, X$, we have
\begin{align*}
\left(\,\sum\limits_{\,i \,\in\, I}\,v_{i}^{\,p}\, \left\|\,\Lambda_{i}\,P_{\,V_{i}}\,U\,P_{\,U\,V_{i}}\,(\,f\,) \,\right\|^{\,p}\,\right)^{1 \,/\, p}& \,=\, \left(\,\sum\limits_{\,i \,\in\, I}\,v_{i}^{\,p}\, \left\|\,\Lambda_{i}\,P_{\,V_{i}}\,U\,(\,f\,)\,\right\|^{\,p}\,\right)^{1 \,/\, p}\\
&\geq\, A\,\left\|\,U\,f\,\right\| \,\geq\, A\,M\,\left\|\,f\,\right\|.
\end{align*}
According to the proof of the Theorem \ref{th1.20}, the upper frame condition is also satisfied.\,This completes the proof.    
\end{proof}

We now give a characterization of a \,$g$-$p$-fusion Bessel sequence in \,$X$.

\begin{theorem}\label{th1.21}
The family \,$\Lambda$\, is a g-p-fusion Bessel sequence in \,$X$\, with respect to \,$\left\{\,X_{i}\,\right\}_{i \,\in\, I}$\, having bound \,$B$\, if and only if the operator given by
\[T \,:\, l^{\,q}\,\left(\,\left\{\,X^{\,\ast}_{i}\,\right\}_{i \,\in\, I}\,\right) \,\to\, X^{\,\ast},\; \;T\,\left(\,\left\{\,g_{i}\,\right\}_{i \,\in\, I}\,\right) \,=\, \sum\limits_{i \,\in\, I}\,v_{\,i}\,P_{\,V_{i}}\,\Lambda^{\,\ast}_{i}\,g_{\,i}\]
is a well-defined, bounded linear operator with \,$\|\,T\,\| \,\leq\, B$.
\end{theorem}

\begin{proof}
First we suppose that \,$\Lambda$\, is a \,$g$-$p$-fusion Bessel sequence in \,$X$\, with respect to \,$\left\{\,X_{i}\,\right\}_{i \,\in\, I}$\, having bound \,$B$.\,Then for any \,$\left\{\,g_{i}\,\right\}_{i \,\in\, I} \,\in\, l^{\,q}\,\left(\,\left\{\,X^{\,\ast}_{i}\,\right\}_{i \,\in\, I}\,\right)$\, and any subset \,$J \,\subset\, I$, we have
\begin{align*}
&\left\|\,\sum\limits_{i \,\in\, J}\,v_{\,i}\,P_{\,V_{i}}\,\Lambda^{\,\ast}_{i}\,g_{\,i}\,\right\| \,=\, \sup\limits_{f \,\in\, X,\, \,\|\,f\,\| \,=\, 1}\,\left|\,\sum\limits_{i \,\in\, J}\,v_{\,i}\,P_{\,V_{i}}\,\Lambda^{\,\ast}_{i}\,g_{\,i}\,(\,f\,)\,\right|\\ 
&\,=\, \sup\limits_{f \,\in\, X,\, \,\|\,f\,\| \,=\, 1}\,\left|\,\sum\limits_{i \,\in\, J}\,g_{\,i}\,v_{\,i}\,\Lambda_{i}\,P_{\,V_{i}}\,(\,f\,)\,\right| \leq\, \sup\limits_{f \,\in\, X,\, \,\|\,f\,\| \,=\, 1}\,\sum\limits_{i \,\in\, J}\,\left\|\,g_{\,i}\,\right\|\,v_{\,i}\,\left\|\,\Lambda_{i}\,P_{\,V_{i}}\,(\,f\,) \,\right\|\\
&\leq\, \sup\limits_{f \,\in\, X,\, \,\|\,f\,\| \,=\, 1}\,\left(\,\sum\limits_{i \,\in\, J}\,\left\|\,g_{\,i}\,\right\|^{\,q}\,\right)^{1 \,/\, q}\,\left(\,\sum\limits_{\,i \,\in\, J}\,v_{j}^{\,p}\, \left\|\,\Lambda_{i}\,P_{\,V_{i}}\,(\,f\,)\,\right\|^{\,p}\,\right)^{1 \,/\, p}\\
&\leq\,B\,\left(\,\sum\limits_{\,i \,\in\, J}\,\left\|\,g_{\,i}\,\right\|^{\,q}\,\right)^{1 \,/\, q}\; \;[\;\text{since $\Lambda$\, is a \,$g$-$p$-fusion Bessel sequence}\;]. 
\end{align*}
This shows that  the series \,$\sum\limits_{i \,\in\, I}\,v_{\,i}\,P_{\,V_{i}}\,\Lambda^{\,\ast}_{i}\,g_{\,i}$\, is unconditionally convergent in \,$X^{\,\ast}$.\,From the above calculation also it follows that
\begin{align*}
&\left\|\,\sum\limits_{i \,\in\, I}\,v_{\,i}\,P_{\,V_{i}}\,\Lambda^{\,\ast}_{i}\,g_{\,i}\,\right\| \,\leq\, B\,\left(\,\sum\limits_{\,i \,\in\, I}\,\left\|\,g_{\,i}\,\right\|^{\,q}\,\right)^{1 \,/\, q}\\
&\Rightarrow\,\left\|\,T\,\left(\,\left\{\,g_{i}\,\right\}_{i \,\in\, I}\,\right)\,\right\| \,\leq\, B\,\left(\,\sum\limits_{\,i \,\in\, I}\,\left\|\,g_{\,i}\,\right\|^{\,q}\,\right)^{1 \,/\, q} \,=\, B\,\left\|\,\left\{\,g_{i}\,\right\}_{i \,\in\, I}\,\right\|_{q}.
\end{align*}
Thus \,$T$\, is bounded and \,$\|\,T\,\| \,\leq\, B$.\\

Conversely, suppose that \,$T$\, is well-defined and bounded linear operator.\,For fixed \,$f \,\in\, X$, consider the mapping \,$F_{f} \,:\, l^{\,q}\,\left(\,\left\{\,X^{\,\ast}_{i}\,\right\}_{i \,\in\, I}\,\right) \,\to\, \mathbb{C}$\, defined by
\[F_{f}\,\left(\,\left\{\,g_{i}\,\right\}_{i \,\in\, I}\,\right) \,=\, T\,\left(\,\left\{\,g_{i}\,\right\}_{i \,\in\, I}\,\right)\,(\,f\,) \,=\, \sum\limits_{i \,\in\, I}\,v_{\,i}\,g_{\,i}\,\Lambda_{i}\,P_{\,V_{i}}\,(\,f\,).\]Then \,$F_{f}$\, is a bounded linear functional on \,$l^{\,q}\,\left(\,\left\{\,X^{\,\ast}_{i}\,\right\}_{i \,\in\, I}\,\right)$, so 
\[\left\{\,v_{\,i}\,\Lambda_{i}\,P_{\,V_{i}}\,(\,f\,)\,\right\} \,\in\, l^{\,p}\,\left(\,\left\{\,X_{i}\,\right\}_{i \,\in\, I}\,\right)\] and 
\[\left\|\,F_{f}\,\left(\,\left\{\,g_{i}\,\right\}_{i \,\in\, I}\,\right)\,\right\| \,\leq\, \|\,T\,\|\,\left\|\,\left\{\,g_{i}\,\right\}_{i \,\in\, I}\,\right\|_{q}\,\|\,f\,\|.\]
Now, by the Hahn-Banach Theorem, there exists \,$\left\{\,g_{i}\,\right\}_{i \,\in\, I} \,\in\, l^{\,q}\,\left(\,\left\{\,X^{\,\ast}_{i}\,\right\}_{i \,\in\, I}\,\right)$\, with \,$\left\|\,\left\{\,g_{i}\,\right\}_{i \,\in\, I}\,\right\|_{q} \,\leq\, 1$\, such that 
\[\left\|\,\left\{\,v_{\,i}\,\Lambda_{i}\,P_{\,V_{i}}\,(\,f\,)\,\right\}\,\right\|_{p} \,=\, \left|\,\sum\limits_{i \,\in\, I}\,v_{\,i}\,g_{\,i}\,\Lambda_{i}\,P_{\,V_{i}}\,(\,f\,)\,\right|.\]Thus
\begin{align*}
&\left(\,\sum\limits_{\,i \,\in\, I}\,v_{i}^{\,p}\, \left\|\,\Lambda_{i}\,P_{\,V_{i}}\,(\,f\,) \,\right\|^{\,p}\,\right)^{1 \,/\, p} \,=\, \left\|\,\left\{\,v_{\,i}\,\Lambda_{j}\,P_{\,V_{i}}\,(\,f\,)\,\right\}\,\right\|_{p}\\
& \,\leq\, \sup\limits_{\left\|\,\left\{\,g_{i}\,\right\}_{i \,\in\, I}\,\right\|_{q} \,\leq\, 1}\,\left|\,\sum\limits_{i \,\in\, I}\,v_{\,i}\,g_{\,i}\,\Lambda_{i}\,P_{\,V_{i}}\,(\,f\,)\,\right| \,=\, \left\|\,F_{f}\,\right\| \,\leq\, \|\,T\,\|\,\|\,f\,\|.
\end{align*}
This completes the proof. 
\end{proof}

\begin{definition}
Let \,$\Lambda$\, be a g-p-fusion frame for \,$X$\, with respect to \,$\left\{\,X_{i}\,\right\}_{i \,\in\, I}$. Then the operator given by
\[U \,:\, X \,\to\, l^{\,p}\,\left(\,\left\{\,X_{i}\,\right\}_{i \,\in\, I}\,\right),\; \;U\,f \,=\, \left\{\,v_{\,i}\,\Lambda_{i}\,P_{\,V_{i}}\,(\,f\,)\,\right\}_{i \,\in\, I}\; \;\forall\; f \,\in\, X.\]
is called the analysis operator and the operator \,$T \,:\, l^{\,q}\,\left(\,\left\{\,X^{\,\ast}_{i}\,\right\}_{i \,\in\, I}\,\right) \,\to\, X^{\,\ast}$, 
\[T\,\left(\,\left\{\,g_{i}\,\right\}_{i \,\in\, I}\,\right) \,=\, \sum\limits_{i \,\in\, I}\,v_{\,i}\,P_{\,V_{i}}\,\Lambda^{\,\ast}_{i}\,g_{\,i}\; \;\forall\, \left\{\,g_{i}\,\right\}_{i \,\in\, I} \,\in\, l^{\,q}\,\left(\,\left\{\,X^{\,\ast}_{i}\,\right\}_{i \,\in\, I}\,\right)\]
is called synthesis operator. 
\end{definition}

\begin{lemma}\label{lm1}
Let \,$\Lambda$\, be a g-p-fusion frame for \,$X$\, with respect to \,$\left\{\,X_{i}\,\right\}_{i \,\in\, I}$.\,Then the analysis operator \,$U$\, has closed range.
\end{lemma}

\begin{proof}
Since \,$\Lambda$\, be a \,$g$-$p$-fusion frame for \,$X$, by the definition of analysis operator \,$U$, the inequality (\ref{eq1}), can be written as \,$A \,\left\|\,f \,\right\| \,\leq\, \left\|\,U\,f \,\right\| \,\leq\, B \,\left\|\,f \,\right\|$.\,Now, it is easy to verify that \,$U$\, is one-to-one, \,$X \,\cong\, \mathcal{R}_{\,U}$\, and hence \,$U$\, has closed range.  
\end{proof}

\begin{lemma}\label{lm1.1}
Let \,$\Lambda$\, be a g-p-fusion frame for \,$X$\, with respect to \,$\left\{\,X_{i}\,\right\}_{i \,\in\, I}$.\,If for each \,$i \,\in\, I,\, \,X_{i}$\, is reflexive then \,$X$\, is reflexive.
\end{lemma}

\begin{proof}
The proof is follows from the lemma \ref{lm1}.
\end{proof}

\begin{theorem}
Let \,$\Lambda$\, be a g-p-fusion Bessel sequence in \,$X$\, with respect to \,$\left\{\,X_{i}\,\right\}_{i \,\in\, I}$.\,Then
\begin{description}
\item[$(i)$]$U^{\,\ast} \,=\, T$.
\item[$(ii)$]If \,$\Lambda$\, has the lower g-p-fusion frame condition and for each \,$i \,\in\, I,\, \,X_{i}$\, is reflexive then \,$T^{\,\ast} \,=\, U$. 
\end{description}
\end{theorem}

\begin{proof}$(i)$
For any \,$f \,\in\, X$\, and \,$\left\{\,g_{i}\,\right\}_{i \,\in\, I} \,\in\, l^{\,q}\,\left(\,\left\{\,X^{\,\ast}_{i}\,\right\}_{i \,\in\, I}\,\right)$, we have
\begin{align*}
&\left<\,U\,f,\, \left\{\,g_{i}\,\right\}_{i \,\in\, I}\right> \,=\, \left<\,\left\{\,v_{\,i}\,\Lambda_{i}\,P_{\,V_{i}}\,(\,f\,)\,\right\}_{i \,\in\, I},\, \left\{\,g_{i}\,\right\}_{i \,\in\, I}\right>\,=\, \sum\limits_{i \,\in\, I}\,\left<\,v_{\,i}\,\Lambda_{i}\,P_{\,V_{i}}\,(\,f\,),\, g_{i}\,\right>,\\
&\left<\,f,\, T\,\left(\,\left\{\,g_{i}\,\right\}_{i \,\in\, I}\,\right)\,\right> \,=\, \left<\,f,\, \sum\limits_{i \,\in\, I}\,v_{\,i}\,P_{\,V_{i}}\,\Lambda^{\,\ast}_{i}\,g_{\,i}\,\right> \,=\, \sum\limits_{i \,\in\, I}\,\left<\,v_{\,i}\,\Lambda_{i}\,P_{\,V_{i}}\,(\,f\,),\, g_{i}\,\right>.
\end{align*}
This shows that \,$U^{\,\ast} \,=\, T$.\\

$(ii)$\; The proof is directly follows from the lemma \ref{lm1}.  
\end{proof}

The following Theorem gives a characterization of a \,$g$-$p$-fusion frame for \,$X$.

\begin{theorem}\label{th1.3}
The family \,$\Lambda$\, is a g-p-fusion frame for \,$X$\, with respect to \,$\left\{\,X_{i}\,\right\}_{i \,\in\, I}$\, if and only if the synthesis operator \,$T$\, is a surjective and bounded linear operator.
\end{theorem}

\begin{proof}
First we consider that \,$\Lambda$\, is a \,$g$-$p$-fusion frame for \,$X$.\,Then by Theorem \ref{th1.21}, \,$T$\, is well-defined and bounded linear operator.\,Since \,$U$\, is one-to-one and \,$U^{\,\ast} \,=\, T$, by Theorem \ref{thm1} \,$T$\, is onto.\\

Conversely, suppose that \,$T$\, is bounded and onto.\,Then by Theorem \ref{th1.21}, \,$\Lambda$\, is a \,$g$-$p$-fusion Bessel sequence in \,$X$.\,Also, by Theorem \ref{thm1}, \,$U$\, has a bounded inverse and this gives the lower \,$g$-$p$-fusion frame condition.\,This completes the proof.  
\end{proof}

We now develop the concept of generalized Riesz basis into the Banach space \,$X$.

\begin{definition}
Let \,$1 \,<\, q \,<\, \infty$.\,The family \,$\Lambda$\, is called a \,$q$-$gf$-Riesz basis for \,$X$\, with respect to \,$\left\{\,X_{i}\,\right\}_{i \,\in\, J}$\, if
\begin{description}
\item[$(i)$]$\Lambda$\, is $gf$-complete, i\,.\,e., \,$\left\{\,f \,:\, \Lambda_{i}\,P_{\,V_{i}}\,(\,f\,) \,=\, 0,\, i \,\in\, I\,\right\} \,=\, \{\,0\,\}$.
\item[$(ii)$]There exist constants \,$0 \,<\, A \,\leq\, B \,<\, \infty$\, such that for any subset \,$J \,\subset\, I$\, and \,$g_{\,i} \,\in\, X^{\,\ast}_{i},\, i \,\in\, J$,
\[A\,\left(\,\sum\limits_{\,i \,\in\, J}\,\left\|\,g_{\,i}\,\right\|^{\,q}\,\right)^{1 \,/\, q} \,\leq\, \left\|\,\sum\limits_{i \,\in\, J}\,v_{\,i}\,P_{\,V_{i}}\,\Lambda^{\,\ast}_{i}\,g_{\,i}\,\right\| \,\leq\, B\,\left(\,\sum\limits_{\,i \,\in\, J}\,\left\|\,g_{\,i}\,\right\|^{\,q}\,\right)^{1 \,/\, q}.\] 
\end{description} 
\end{definition}

Next theorem establish a relationship between \,$q$-$gf$-Riesz basis and the synthesis operator \,$T$.

\begin{theorem}\label{th1.4}
The family \,$\Lambda$\, is a \,$q$-$gf$-Riesz basis for \,$X$\, with respect to \,$\left\{\,X_{i}\,\right\}_{i \,\in\, I}$\, having bounds \,$A$\, and \,$B$\, if and only if the synthesis operator \,$T$\, is a bounded linear and invertible such that
\begin{equation}\label{eqn1.3}
A\,\|\,g\,\| \,\leq\, \|\,T\,g\,\| \,\leq\, B\,\|\,g\,\|
\end{equation}
for any \,$g \,=\, \,\left\{\,g_{i}\,\right\}_{i \,\in\, I} \,\in\, l^{\,q}\,\left(\,\left\{\,X^{\,\ast}_{i}\,\right\}_{i \,\in\, I}\,\right)$.
\end{theorem}

\begin{proof}
Suppose \,$\Lambda$\, is a \,$q$-$gf$-Riesz basis for \,$X$\, with respect to \,$\left\{\,X_{i}\,\right\}_{i \,\in\, J}$\, having bounds \,$A$\, and \,$B$.\,Then from the definition of \,$q$-$gf$-Riesz basis, it is easy to verify that \,$\sum\limits_{i \,\in\, I}\,v_{\,i}\,P_{\,V_{i}}\,\Lambda^{\,\ast}_{i}\,g_{\,i}$\, converges unconditionally for all \,$\left\{\,g_{i}\,\right\}_{i \,\in\, I} \,\in\, l^{\,q}\,\left(\,\left\{\,X^{\,\ast}_{i}\,\right\}_{i \,\in\, I}\,\right)$,
\[A\,\left(\,\sum\limits_{\,i \,\in\, I}\,\left\|\,g_{\,i}\,\right\|^{\,q}\,\right)^{1 \,/\, q} \,\leq\, \left\|\,\sum\limits_{i \,\in\, I}\,v_{\,i}\,P_{\,V_{i}}\,\Lambda^{\,\ast}_{i}\,g_{\,i}\,\right\| \,\leq\, B\,\left(\,\sum\limits_{\,i \,\in\, I}\,\left\|\,g_{\,i}\,\right\|^{\,q}\,\right)^{1 \,/\, q}\]
and this implies that \,$T$\, is bounded, one-to-one and \,$A\,\|\,g\,\| \,\leq\, \|\,T\,g\,\| \,\leq\, B\,\|\,g\,\|$.\\

Conversely, suppose that the operator \,$T$\, is a bounded linear and invertible  operator from \,$l^{\,q}\,\left(\,\left\{\,X^{\,\ast}_{i}\,\right\}_{i \,\in\, I}\,\right)$ \, onto \,$X^{\,\ast}$\, and satisfying (\ref{eqn1.3}).\,Then by Theorem \ref{th1.3}, \,$\Lambda$\, is a $g$-$p$-fusion frame for \,$X$\, with respect to \,$\left\{\,X_{i}\,\right\}_{i \,\in\, J}$\, having bounds \,$A$\, and \,$B$.\,Now, for \,$f \,\in\, \left\{\,f \,:\, \Lambda_{i}\,P_{\,V_{i}}\,(\,f\,) \,=\, 0,\, i \,\in\, I\,\right\}$, we have
\[A\,\|\,f\,\| \,\leq\, \left(\,\sum\limits_{\,i \,\in\, I}\,v_{i}^{\,p}\, \left\|\,\Lambda_{i}\,P_{\,V_{i}}\,(\,f\,) \,\right\|^{\,p}\,\right)^{1 \,/\, p} \,=\, 0\; \;\Rightarrow\, f \,=\, 0.\]
Therefore, we obtain that \,$\left\{\,f \,:\, \Lambda_{i}\,P_{\,V_{i}}\,(\,f\,) \,=\, 0,\, i \,\in\, I\,\right\} \,=\, \{\,0\,\}$.\,Hence, \,$\Lambda$\, is a \,$q$-$gf$-Riesz basis for \,$X$\, with respect to \,$\left\{\,X_{i}\,\right\}_{i \,\in\, I}$\, having bounds \,$A$\, and \,$B$.           
\end{proof}

\begin{remark}
Let \,$\Lambda$\, be a \,$q$-$gf$-Riesz basis for \,$X$\, with respect to \,$\left\{\,X_{i}\,\right\}_{i \,\in\, I}$\, having bounds \,$A$\, and \,$B$.\,Then \,$\Lambda$\, is also a g-p-fusion frame for \,$X$\, with respect to \,$\left\{\,X_{i}\,\right\}_{i \,\in\, I}$\, having same bounds.
\end{remark}

\begin{proof}
By Theorem \ref{th1.4}, \,$T$\, is a bounded linear invertible operator with \,$\|\,T\,\| \,\leq\, B$\, and \,$\left\|\,T^{\,-\, 1}\,\right\| \,\leq\, A^{\,-\, 1}$.\,It is easy to verify that \,$\left\|\,(\,T^{\,\ast}\,)^{\,-\, 1}\,\right\|^{\,-\, 1} \,\geq\, A$.\,Then by Theorem \ref{th1.3}, \,$\Lambda$\, is a \,$g$-$p$-fusion frame for \,$X$\, with respect to \,$\left\{\,X_{i}\,\right\}_{i \,\in\, I}$\, having bounds \,$A$\, and \,$B$. 
\end{proof}

\begin{theorem}
Let \,$\left\{\,X_{i}\,\right\}_{i \,\in\, I}$\, be a sequence of reflexive Banach spaces and \,$\Lambda$\, be a g-p-fusion frame for \,$X$\, with respect to \,$\left\{\,X_{i}\,\right\}_{i \,\in\, I}$.\,Then the following are equivalent
\begin{description}
\item[$(i)$]\,$\Lambda$\, is a \,$q$-$gf$-Riesz basis for \,$X$\, with respect to \,$\left\{\,X_{i}\,\right\}_{i \,\in\, I}$.
\item[$(ii)$]If for any \,$g \,=\, \,\left\{\,g_{i}\,\right\}_{i \,\in\, I} \,\in\,  l^{\,q}\,\left(\,\left\{\,X^{\,\ast}_{i}\,\right\}_{i \,\in\, I}\,\right),\, \, \sum\limits_{i \,\in\, I}\,v_{\,i}\,P_{\,V_{i}}\,\Lambda^{\,\ast}_{i}\,g_{\,i} \,=\, 0$, then \,$g_{\,i} \,=\, 0\; \,\forall\, i \,\in\, I$. 
\item[$(iii)$]$\mathcal{R}\,(\,U\,) \,=\, l^{\,p}\,\left(\,\left\{\,X_{i}\,\right\}_{i \,\in\, I}\,\right)$
\end{description} 
\end{theorem}

\begin{proof}
From the definition of \,$q$-$gf$-Riesz basis, it is easy to verify \,$(i)  \,\Rightarrow\, (ii)$.\\

$(ii)  \,\Rightarrow\, (i)$ Suppose that \,$(ii)$\, holds.\,Since \,$\Lambda$\, be a \,$g$-$p$-fusion frame for \,$X$\, with respect to \,$\left\{\,X_{i}\,\right\}_{i \,\in\, I}$, by Theorem \ref{th1.3}, the operator \,$T$\, is bounded linear and surjective.\,Also by condition \,$(ii)$, it is easy to verify that \,$T$\, is injective.\,Hence, by Theorem \ref{th1.4}, \,$\Lambda$\, is a \,$q$-$gf$-Riesz basis for \,$X$\, with respect to \,$\left\{\,X_{i}\,\right\}_{i \,\in\, I}$.\\

$(i) \,\Rightarrow\, (iii)$\, and \,$(iii) \,\Rightarrow\, (i)$\, are directly follows from the Theorem 3.16 of \cite{M}. 
\end{proof}

\section{Perturbation of $g$-$p$-fusion frame}

\smallskip\hspace{.6 cm}In this section, the stability of \,$g$-$p$-fusion frame in \,$X$\, is presented.

\begin{theorem}
Let \,$\Lambda$\, be a g-p-fusion frame for \,$X$\, with respect to \,$\left\{\,X_{i}\,\right\}_{i \,\in\, I}$\, having bounds \,$A$\, and \,$B$.\,Suppose that \,$\Gamma_{i} \,\in\, \mathcal{B}\,\left(\,X,\,X_{i}\,\right),\, i \,\in\, I$\, such that
\[\left(\sum\limits_{\,i \,\in\, I}v_{i}^{\,p}\left\|\,\left(\,\Lambda_{i}\,P_{\,V_{i}} \,-\, \Gamma_{i}\,P_{\,W_{i}}\right)\,(\,f\,)\,\right\|^{\,p}\right)^{1 \,/\, p} \leq \lambda_{1}\,\left(\sum\limits_{\,i \,\in\, I}\,v_{i}^{\,p}\, \left\|\,\Lambda_{i}\,P_{\,V_{i}}\,(\,f\,) \,\right\|^{\,p}\right)^{1 \,/\, p} +\]
\begin{equation}\label{eq1.3}
+\, \lambda_{2}\left(\,\sum\limits_{\,i \,\in\, I}\,v_{i}^{\,p}\, \left\|\,\Gamma_{i}\,P_{\,W_{i}}\,(\,f\,) \,\right\|^{\,p}\,\right)^{1 \,/\, p} \,+\, \mu\,\|\,f\,\|\; \;\forall\, f \,\in\, X,
\end{equation}
where \,$\lambda_{1},\,\lambda_{2} \,\in\, (\,-\,1,\, 1\,)$\, and \,$-\,\left(\,1 \,+\, \lambda_{1}\,\right)\,B \,\leq\, \mu \,\leq\, \left(\,1 \,-\, \lambda_{1}\,\right)\,A$.\,Then \,$\Gamma \,=\, \left\{\,\left(\,W_{i},\, \Gamma_{i},\, v_{i}\,\right)\,\right\}_{i \,\in\, I}$\, is a g-p-fusion frame for \,$X$\, with respect to \,$\left\{\,X_{i}\,\right\}_{i \,\in\, I}$.  
\end{theorem}

\begin{proof}
For each \,$f \,\in\, X$, by Minkowski inequality, we have
\begin{align*}
&\left(\,\sum\limits_{\,i \,\in\, I}\,v_{i}^{\,p}\, \left\|\,\Gamma_{i}\,P_{\,W_{i}}\,(\,f\,) \,\right\|^{\,p}\,\right)^{1 \,/\, p} \,\leq\, \left(\,\sum\limits_{\,i \,\in\, I}\,v_{i}^{\,p}\,\left\|\,\left(\,\Lambda_{i}\,P_{\,V_{i}} \,-\, \Gamma_{i}\,P_{\,W_{i}}\,\right)\,(\,f\,)\,\right\|^{\,p}\,\right)^{1 \,/\, p} +\\
&\hspace{4cm}+\,\left(\,\sum\limits_{\,i \,\in\, I}\,v_{i}^{\,p}\, \left\|\,\Lambda_{i}\,P_{\,V_{i}}\,(\,f\,) \,\right\|^{\,p}\,\right)^{1 \,/\, p}\\
&\leq\left(\,1 \,+\, \lambda_{1}\,\right)\left(\sum\limits_{\,i \,\in\, I}v_{i}^{\,p}\left\|\,\Lambda_{i}\,P_{\,V_{i}}\,(\,f\,) \,\right\|^{\,p}\right)^{1 \,/\, p} \,+\, \lambda_{2}\left(\sum\limits_{\,i \,\in\, I}v_{i}^{\,p}\, \left\|\,\Gamma_{i}\,P_{\,W_{i}}\,(\,f\,) \,\right\|^{\,p}\right)^{1 \,/\, p} +\, \mu\,\|\,f\,\|
\end{align*}
Therefore, since \,$\Lambda$\, is a \,$g$-$p$-fusion frame, we have
\begin{align*}
&\left(\sum\limits_{\,i \,\in\, I}v_{i}^{\,p}\left\|\,\Gamma_{i}\,P_{\,W_{i}}\,(\,f\,) \,\right\|^{\,p}\right)^{1 \,/\, p}\\
& \leq \dfrac{1 \,+\, \lambda_{1}}{1 \,-\, \lambda_{2}}\,\left(\sum\limits_{\,i \,\in\, I}v_{i}^{\,p}\left\|\,\Lambda_{i}\,P_{\,V_{i}}\,(\,f\,) \,\right\|^{\,p}\right)^{1 \,/\, p} +\, \dfrac{\mu}{\left(\,1 \,-\, \lambda_{2}\,\right)}\,\|\,f\,\|\\
&\leq\, \left[\,\left(\,\dfrac{1 \,+\, \lambda_{1}}{1 \,-\, \lambda_{2}}\,\right)\,B \,+\, \dfrac{\mu}{\left(\,1 \,-\, \lambda_{2}\,\right)}\,\right]\,\|\,f\,\| \,=\, \left[\,\dfrac{B\,\left(\,1 \,+\, \lambda_{1}\,\right) \,+\, \mu}{1 \,-\, \lambda_{2}}\,\right]\,\|\,f\,\|.
\end{align*} 
On the other hand,
\begin{align*}
&\left(\,\sum\limits_{\,i \,\in\, I}\,v_{i}^{\,p}\, \left\|\,\Lambda_{i}\,P_{\,V_{i}}\,(\,f\,) \,\right\|^{\,p}\,\right)^{1 \,/\, p} \,-\, \left(\,\sum\limits_{\,i \,\in\, I}\,v_{i}^{\,p}\, \left\|\,\Gamma_{i}\,P_{\,W_{i}}\,(\,f\,) \,\right\|^{\,p}\,\right)^{1 \,/\, p}\\
&\hspace{2cm}\leq\, \left(\,\sum\limits_{\,i \,\in\, I}\,v_{i}^{\,p}\,\left\|\,\left(\,\Lambda_{i}\,P_{\,V_{i}} \,-\, \Gamma_{i}\,P_{\,W_{i}}\,\right)\,(\,f\,)\,\right\|^{\,p}\,\right)^{1 \,/\, p}.
\end{align*}
Now using (\ref{eq1.3}), we obtain
\begin{align*}
&\left(\,1 \,+\, \lambda_{2}\,\right)\left(\sum\limits_{\,i \,\in\, I}\,v_{i}^{\,p}\, \left\|\,\Gamma_{i}\,P_{\,W_{i}}\,(\,f\,)\,\right\|^{\,p}\right)^{1 \,/\, p}\\
& \geq \left(\,1 \,-\, \lambda_{1}\,\right)\left(\sum\limits_{\,i \,\in\, I}\,v_{i}^{\,p}\left\|\,\Lambda_{i}\,P_{\,V_{i}}\,(\,f\,) \,\right\|^{\,p}\right)^{1 \,/\, p} -\, \mu\,\|\,f\,\|\\
&\geq\, \left[\,\left(\,1 \,-\, \lambda_{1}\,\right)\,A \,-\, \mu\,\right]\,\|\,f\,\|\; \;[\;\text{since $\Lambda$ is a $g$-$p$-fusion frame}\;]\\
&\Rightarrow\, \left(\,\sum\limits_{\,i \,\in\, I}\,v_{i}^{\,p}\, \left\|\,\Gamma_{i}\,P_{\,W_{i}}\,(\,f\,) \,\right\|^{\,p}\,\right)^{1 \,/\, p} \,\geq\, \left[\,\dfrac{A\,\left(\,1 \,-\, \lambda_{1}\,\right) \,-\, \mu}{1 \,+\, \lambda_{2}}\,\right]\,\|\,f\,\|.
\end{align*}
Hence, \,$\Gamma$\, is a \,$g$-$p$-fusion frame for \,$X$\, with respect to \,$\left\{\,X_{i}\,\right\}_{i \,\in\, I}$\, having bounds
 \[\left[\,\dfrac{A\,\left(\,1 \,-\, \lambda_{1}\,\right) \,-\, \mu}{1 \,+\, \lambda_{2}}\,\right]\, \,\text{and}\, \,\left[\,\left(\,\dfrac{1 \,+\, \lambda_{1}}{1 \,-\, \lambda_{2}}\,\right)\,B \,+\, \dfrac{\mu}{\left(\,1 \,-\, \lambda_{2}\,\right)}\,\right].\]\\This completes the proof.  
\end{proof}

\begin{theorem}
Let \,$\Lambda$\, be a g-p-fusion frame for \,$X$\, with respect to \,$\left\{\,X_{i}\,\right\}_{i \,\in\, I}$\, having bounds \,$A$\, and \,$B$.\,Suppose that \,$\Gamma_{i} \,\in\, \mathcal{B}\,\left(\,X,\,X_{i}\,\right),\, i \,\in\, I$\, such that
\[\left(\,\sum\limits_{\,i \,\in\, I}\,v_{i}^{\,p}\,\left\|\,\left(\,\Lambda_{i}\,P_{\,V_{i}} \,-\, \Gamma_{i}\,P_{\,W_{i}}\,\right)\,(\,f\,)\,\right\|^{\,p}\,\right)^{1 \,/\, p} \,\leq\, R\,\|\,f\,\|\; \;\forall\; f \,\in\, X.\]
where \,$0 \,<\, R \,<\, A$.\,Then \,$\Gamma \,=\, \left\{\,\left(\,W_{i},\, \Gamma_{i},\, v_{i}\,\right)\,\right\}_{i \,\in\, I}$\, is a g-p-fusion frame for \,$X$\, with respect to \,$\left\{\,X_{i}\,\right\}_{i \,\in\, I}$\, having bounds \,$(\,A \,-\, R\,)$\, and \,$(\,B \,+\, R\,)$.
\end{theorem}

\begin{proof}
By Minkowski inequality, for \,$f \,\in\, X$, we get
\begin{align*}
\left(\sum\limits_{\,i \,\in\, I}v_{i}^{\,p}\left\|\,\Gamma_{i}\,P_{\,W_{i}}\,(\,f\,) \,\right\|^{\,p}\,\right)^{1 \,/\, p}& \leq \left(\sum\limits_{\,i \,\in\, I}\,v_{i}^{\,p}\,\left\|\,\left(\,\Lambda_{i}\,P_{\,V_{i}} \,-\, \Gamma_{i}\,P_{\,W_{i}}\,\right)\,(\,f\,)\,\right\|^{\,p}\right)^{1 \,/\, p} \,+\\
&+\,\left(\,\sum\limits_{\,i \,\in\, I}\,v_{i}^{\,p}\, \left\|\,\Lambda_{i}\,P_{\,V_{i}}\,(\,f\,) \,\right\|^{\,p}\,\right)^{1 \,/\, p}\\
&\leq\, (\,B \,+\, R\,)\,\|\,f\,\|\; \;[\;\text{since $\Lambda$ is $g$-$p$-fusion frame}\;].
\end{align*}
On the other hand,
\begin{align*}
\left(\,\sum\limits_{\,i \,\in\, I}\,v_{i}^{\,p}\, \left\|\,\Gamma_{i}\,P_{\,W_{i}}\,(\,f\,) \,\right\|^{\,p}\,\right)^{1 \,/\, p}& \,\geq\, \left(\,\sum\limits_{\,i \,\in\, I}\,v_{i}^{\,p}\, \left\|\,\Lambda_{i}\,P_{\,V_{i}}\,(\,f\,) \,\right\|^{\,p}\,\right)^{1 \,/\, p} \,-\\
&-\, \left(\,\sum\limits_{\,i \,\in\, I}\,v_{i}^{\,p}\,\left\|\,\left(\,\Lambda_{i}\,P_{\,V_{i}} \,-\, \Gamma_{i}\,P_{\,W_{i}}\,\right)\,(\,f\,)\,\right\|^{\,p}\,\right)^{1 \,/\, p}\\
& \,\geq\, (\,A \,-\, R\,)\|\,f\,\|.
\end{align*}
Hence, \,$\Gamma$\, is a $g$-$p$-fusion frame for \,$X$\, with bounds \,$(\,A \,-\, R\,)$\, and \,$(\,B \,+\, R\,)$. \\This completes the proof.   
\end{proof}

We end this section by constructing \,$g$-$p$-fusion frames in Cartesian product of Banach spaces and tensor product of Banach spaces.\\

Let \,$\left(\,X,\, \left\|\,\cdot\,\right\|_{X}\,\right)$\, and \,$\left(\,Y,\, \left\|\,\cdot\,\right\|_{Y}\,\right)$\, be two Banach spaces.\,Then the Cartesian product of \,$X$\, and \,$Y$\, is denoted by \,$X \,\oplus\, Y$\, and defined to be an Banach space with respect to the norm
\begin{equation}\label{eqp1}
\left\|\,f \,\oplus\, g\,\right\|^{\,p} \,=\, \|\,f\,\|^{\,p}_{X} \,+\, \|\,g\,\|^{\,p}_{Y},
\end{equation}
for all \,$f \,\in\, X\; \;\text{and}\; \,g \,\in\, Y$.\,Now, if \,$U \,\in\, \mathcal{B}\,(\,X,\,X_{i}\,)$\, and \,$V \,\in\, \mathcal{B}\,(\,Y,\,Y_{i}\,)$, then for all \,$f \,\in\, X$\, and \,$g \,\in\, Y$, we define
\[U \,\oplus\, V \,\in\, \mathcal{B}\left(\,X \,\oplus\, Y,\, X_{i}\,\oplus\, Y_{i}\,\right)\; \;\text{by}\; \;(\,U \,\oplus\, V\,)\,(\,f \,\oplus\, g\,) \,=\, U\,f \,\oplus\, V\,g,\]
\[P_{V_{i} \,\oplus\, W_{i}}\,(\,f \,\oplus\, g\,) \,=\, P_{\,V_{i}}\,f \,\oplus\, P_{\,W_{i}}\,g,\]
where \,$\left\{\,Y_{i}\,\right\}_{i \,\in\, I}$\, is a another sequence of Banach spaces and \,$\left\{\,W_{i}\,\right\}_{i \,\in\, I}$\, is the collection of closed subspaces of \,$Y$ and \,$P_{\,W_{i}}$\, are the linear projections of \,$Y$\, onto \,$W_{i}$\, such that \,$P_{\,W_{i}}\,(\,X\,) \,=\, W_{i}$, for \,$i \,\in\, I$.   

\begin{theorem}
Let \,$\Lambda$\, be a g-p-fusion frame for \,$X$\, with respect to \,$\left\{\,X_{i}\,\right\}_{i \,\in\, I}$\, having bounds \,$A,\, B$\, and \,$\Gamma \,=\, \left\{\,\left(\,W_{i},\, \Gamma_{i},\, v_{i}\,\right)\,\right\}_{i \,\in\, I}$\, be a g-p-fusion frame for \,$Y$\, with respect to \,$\left\{\,Y_{i}\,\right\}_{i \,\in\, I}$\, having bounds \,$C,\,D$, where \,$\Gamma_{i} \,\in\, \mathcal{B}\left(\,Y,\, Y_{i}\,\right)$\, for each \,$i \,\in\, I$.\,Then \,$\Lambda \,\oplus\, \Gamma \,=\, \left\{\,\left(\,V_{i} \,\oplus\, W_{i},\, \Lambda_{i} \,\oplus\, \Gamma_{i},\, v_{i}\,\right)\,\right\}_{i \,\in\, I}$\, is a g-p-fusion frame for \,$X \,\oplus\, Y$\, with respect to \,$\left\{\,X_{i} \,\oplus\, Y_{i}\,\right\}_{i \,\in\, I}$\, having bounds \,$\min\left(\,A^{\,p},\,C^{\,p}\,\right)$\, and \,$\max\left(\,B^{\,p},\, D^{\,p}\,\right)$.    
\end{theorem}

\begin{proof}
Since \,$\Lambda$\, and \,$\Gamma$\, are \,$g$-$p$-fusion frames for \,$X$\, and \,$Y$, respectively, 
\begin{equation}\label{eqp1.1}
A^{\,p} \,\left \|\,f \,\right \|_{X}^{\,p} \,\leq\, \sum\limits_{\,i \,\in\, I}\,v_{i}^{\,p}\,\left\|\,\Lambda_{i}\,P_{\,V_{i}}\,(\,f\,) \,\right\|_{X}^{\,p} \,\leq\, B^{\,p}\,\left\|\,f\, \right\|_{X}^{\,p}\, \;\forall\, f \,\in\, X
\end{equation}
\begin{equation}\label{eqp1.2}
C^{\,p}\,\left \|\,g \,\right \|_{Y}^{\,p} \,\leq\, \sum\limits_{\,i \,\in\, I}\,v_{i}^{\,p}\,\left\|\,\Gamma_{i}\,P_{\,W_{i}}\,(\,g\,) \,\right\|_{Y}^{\,p} \,\leq\, D^{\,p}\,\left\|\,g\, \right\|_{Y}^{\,p}\, \,\forall\, g \,\in\, Y.
\end{equation}
Adding (\ref{eqp1.1}) and (\ref{eqp1.2}) and then using (\ref{eqp1}), we get
\begin{align*}
&A^{\,p}\,\left \|\,f \,\right \|_{X}^{\,p} \,+\, C \,\left \|\,g \,\right \|_{Y}^{\,p} \,\leq\, \sum\limits_{\,i \,\in\, I}\,v_{i}^{\,p}\,\left\|\,\Lambda_{j}\,P_{\,V_{i}}\,(\,f\,) \,\right\|_{X}^{\,p} \,+\, \sum\limits_{\,i \,\in\, I}\,v_{i}^{\,p}\,\left\|\,\Gamma_{i}\,P_{\,W_{i}}\,(\,g\,) \,\right\|_{Y}^{\,p}\\
&\hspace{3.7cm}\,\leq\, B^{\,p}\,\left\|\,f\, \right\|_{X}^{\,p} \,+\, D^{\,p}\,\left\|\,g\, \right\|_{Y}^{\,p}.\\
&\Rightarrow \min\left(\,A^{\,p},\,C^{\,p}\,\right)\left\{\left \|\,f \,\right \|_{X}^{\,p} \,+\, \,\left \|\,g \,\right \|_{Y}^{\,p}\right\} \leq \sum\limits_{\,i \,\in\, I}v_{i}^{\,p}\,\left(\left\|\,\Lambda_{i}\,P_{\,V_{i}}\,(\,f\,) \,\right\|_{X}^{\,p} \,+\, \left\|\,\Gamma_{i}\,P_{\,W_{i}}\,(\,g\,) \,\right\|_{Y}^{\,p}\right)\\
&\hspace{3.7cm}\leq\, \max\left(\,B^{\,p},\, D^{\,p}\,\right)\,\left\{\,\left \|\,f \,\right \|_{X}^{\,p} \,+\, \,\left \|\,g \,\right \|_{Y}^{\,p}\,\right\}.\\
&\Rightarrow\, \min\left(\,A^{\,p},\,C^{\,p}\,\right)\,\left\|\,f \,\oplus\, g\,\right\|^{\,p} \,\leq\, \sum\limits_{\,i \,\in\, I}\,v_{i}^{\,p}\,\left\|\,\Lambda_{i}\,P_{\,V_{i}}\,(\,f\,) \,\oplus\, \Gamma_{i}\,P_{\,W_{i}}\,(\,g\,)\,\right\|^{\,p}\\
&\hspace{3.7cm}\,\leq\, \max\left(\,B^{\,p},\, D^{\,p}\,\right)\,\left\|\,f \,\oplus\, g\,\right\|^{\,p}.\\
&\Rightarrow\, \min\left(\,A^{\,p},\,C^{\,p}\,\right)\,\left\|\,f \,\oplus\, g\,\right\|^{\,p} \,\leq\, \sum\limits_{\,i \,\in\, I}\,v_{i}^{\,p}\,\left\|\,\left(\,\Lambda_{i} \,\oplus\, \Gamma_{i}\,\right)\,\left(\,P_{\,V_{i}} \,\oplus\, P_{\,W_{i}}\,\right)\,(\,f \,\oplus\, g\,)\,\right\|^{\,p}\\
&\hspace{3.7cm} \,\leq\, \max\left(\,B^{\,p},\, D^{\,p}\,\right)\,\left\|\,f \,\oplus\, g\,\right\|^{\,p}\; \;\forall\; f \,\oplus\, g \,\in\, X \,\oplus\, Y.\\
&\Rightarrow\, \min\left(\,A^{\,p},\,C^{\,p}\,\right)\,\left\|\,f \,\oplus\, g\,\right\|^{\,p} \,\leq\, \sum\limits_{\,i \,\in\, I}\,v_{i}^{\,p}\,\left\|\,\left(\,\Lambda_{i} \,\oplus\, \Gamma_{i}\,\right)\,P_{\,V_{i} \,\oplus\, W_{i}}\,(\,f \,\oplus\, g\,)\,\right\|^{\,p}\\
&\hspace{3.7cm} \,\leq\, \max\left(\,B^{\,p},\, D^{\,p}\,\right)\,\left\|\,f \,\oplus\, g\,\right\|^{\,p}\; \;\forall\; f \,\oplus\, g \,\in\, X \,\oplus\, Y.
\end{align*}
Thus, \,$\Lambda \,\oplus\, \Gamma$\, is a \,$g$-$p$-fusion frame for \,$X \,\oplus\, Y$\, with respect to \,$\left\{\,X_{i} \,\oplus\, Y_{i}\,\right\}_{i \,\in\, I}$\, having bounds \,$\min\left(\,A^{\,p},\,C^{\,p}\,\right)$\, and \,$\max\left(\,B^{\,p},\, D^{\,p}\,\right)$.\,This completes the proof.    
\end{proof}

The tensor product of \,$X$\, and \,$Y$\, is denoted by \,$X \,\otimes\, Y$\, and it is defined to be an normed space with respect to the norm 
\begin{equation}\label{eqp1.3}
\left\|\,f \,\otimes\, g\,\right\|^{\,p} \,=\, \|\,f\,\|^{\,p}_{X}\;\|\,g\,\|^{\,p}_{Y},
\end{equation}
for all \,$f \,\in\, X\; \;\text{and}\; \,g \,\in\, Y$.\,Then it is easy to verify that \,$X \,\otimes\, Y$\, is complete with respect to the above norm.\,Therefore, \,$X \,\otimes\, Y$\, is a Banach space. 

\begin{remark}
Let \,$U,\, U^{\,\prime} \,\in\, \mathcal{B}\,(\,X,\,X_{i}\,)$\, and \,$V,\, V^{\,\prime} \,\in\, \mathcal{B}\,(\,Y,\,Y_{j}\,)$, for \,$i \,\in\, I$\, and \,$j \,\in\, J$.\,Then for \,$U \,\otimes\, V,\, U^{\,\prime} \,\otimes\, V^{\,\prime} \,\in\, \mathcal{B}\left(\,X \,\otimes\, Y,\, X_{i}\,\otimes\, Y_{j}\,\right)$, we define
\begin{description}
\item[$(i)$]$\left(\,U \,\otimes\, V\,\right)\,(\,f \,\otimes\, g\,) \,=\, U\,f \,\otimes\, V\,g\; \;\text{for all}\; \;f \,\in X,\; g \,\in\, Y$.
\item[$(ii)$]$\left(\,U \,\otimes\, V\,\right)\,\left(\,U^{\,\prime} \,\otimes\, V^{\,\prime}\,\right) \,=\, U\,U^{\,\prime} \,\otimes\, V\,V^{\,\prime}$.
\item[$(iii)$]$P_{V_{i} \,\otimes\, W_{j}}\,(\,f \,\otimes\, g\,) \,=\, P_{\,V_{i}}\,f \,\otimes\, P_{\,W_{j}}\,g\; \;\text{for all}\; \;f \,\in X,\; g \,\in\, Y$.
\end{description}
\end{remark}

\begin{remark}
Let \,$\left\{\,v_{\,i}\,\right\}_{\, i \,\in\, I},\;\left\{\,w_{\,j}\,\right\}_{ j \,\in\, J}$\, be two families of positive weights i\,.\,e., \,$v_{\,i} \,>\, 0\, \;\forall\; i \,\in\, I,\; \,w_{\,j} \,>\, 0\, \;\forall\; j \,\in\, J$\, and \,$\Lambda_{i} \,\otimes\, \Gamma_{j} \,\in\, \mathcal{B}\,(\,X \,\otimes\, Y,\, X_{i} \,\otimes\, Y_{j}\,)$\, for each \,$i \,\in\, I$\, and \,$j \,\in\, J$.\;Then according to the definition (\ref{deff1}), the family \,$\Lambda \,\otimes\, \Gamma \,=\, \left\{\,\left(\,V_{i} \,\otimes\, W_{j},\, \Lambda_{i} \,\otimes\, \Gamma_{j},\, v_{i}\,w_{j}\,\right)\,\right\}_{\,i,\,j}$\, is said to be a g-p-fusion frame for \,$X \,\otimes\, Y$\, with respect to \,$\left\{\,X_{i} \,\otimes\, Y_{j}\,\right\}_{\,i,\,j}$\, if there exist constants \,$A,\,B \,>\, 0$\, such that 
\[A\left\|\,f \,\otimes\, g\,\right\| \leq \left(\sum\limits_{i,\, j}v^{\,p}_{\,i}\,w^{\,p}_{\,j}\,\left\|\left(\,\Lambda_{i} \,\otimes\, \Gamma_{j}\,\right)\,P_{\,V_{\,i} \,\otimes\, W_{\,j}}\,(\,f \,\otimes\, g\,)\,\right\|^{\,P}\right)^{1 \,/\, p} \leq B \left\|\,f \,\otimes\, g\,\right\|\]
for all \,$f \,\otimes\, g \,\in\, X \,\otimes\, Y$.\;The constants \,$A$\, and \,$B$\, are called the frame bounds.
\end{remark}

\begin{theorem}\label{th1.2}
The family \,$\Lambda \,\otimes\, \Gamma$\, is a g-p-fusion frame for \,$X \,\otimes\, Y$\, with respect to \,$\left\{\,X_{i} \,\otimes\, Y_{j}\,\right\}_{\,i,\,j}$\, if and only if \,$\Lambda$\, is a g-p-fusion frames for \,$X$\,  with respect to \,$\left\{\,X_{i}\,\right\}_{\, i \,\in\, I}$\, and \,$\Gamma$\, is a g-p-fusion frames for \,$Y$\, with respect to \,$\left\{\,Y_{j}\,\right\}_{\, j \,\in\, J}$. 
\end{theorem}

\begin{proof}
First we suppose that \,$\Lambda \,\otimes\, \Gamma$\, is a \,$g$-$p$-fusion frame for \,$H \,\otimes\, K$\, with respect to \,$\left\{\,H_{i} \,\otimes\, K_{j}\,\right\}_{\,i,\,j}$.\;Then there exist constants \,$A,\, B \,>\, 0$\, such that for all \,$f \,\otimes\, g \,\in\, H \,\otimes\, K \,-\, \{\,\theta \,\otimes\, \theta\,\}$, we have
\begin{align*}
&A\left\|\,f \,\otimes\, g\,\right\| \leq \left(\sum\limits_{i,\, j}\,v^{\,p}_{\,i}\,w^{\,p}_{\,j}\,\left\|\left(\,\Lambda_{i} \,\otimes\, \Gamma_{j}\,\right)\,P_{\,V_{\,i} \,\otimes\, W_{\,j}}\,(\,f \,\otimes\, g\,)\,\right\|^{\,P}\right)^{1 \,/\, p} \leq B\left\|\,f \,\otimes\, g\,\right\|\\
&\Rightarrow A\left\|\,f \,\otimes\, g\,\right\| \leq \left(\sum\limits_{i,\, j}v^{\,p}_{\,i}\,w^{\,p}_{\,j}\left\|\,\Lambda_{i}\,P_{\,V_{i}}\,(\,f\,) \otimes\, \Gamma_{j}\,P_{\,W_{j}}\,(\,g\,) \,\right\|^{\,p}\right)^{1 \,/\, p} \leq B\left\|\,f \,\otimes\, g\,\right\|.\\
&\Rightarrow A\left\|\,f\,\right\|_{X}\left\|\,g\,\right\|_{Y} \,\leq\, \left(\sum\limits_{\,i \,\in\, I} v_{i}^{\,p}\left\|\,\Lambda_{i}\,P_{\,V_{i}}\,(\,f\,) \,\right\|_{X}^{\,p}\right)^{1 \,/\, p}\left(\sum\limits_{\,j \,\in\, J} w_{j}^{\,p} \left\|\,\Gamma_{j}\,P_{\,W_{j}}\,(\,g\,) \,\right\|_{Y}^{\,p}\right)^{1 \,/\, p}\\
&\hspace{3cm} \leq\, B\,\left\|\,f\,\right\|_{X}\,\left\|\,g\,\right\|_{Y}\; \;[\,\text{by (\ref{eqp1.3})}\,].
\end{align*}
Since \,$f \,\otimes\, g$\, is non-zero vector, \,$f$\, and \,$g$\, are also non-zero vectors and therefore \,$\sum\limits_{\,i \,\in\, I}\, v_{i}^{\,p} \,\left\|\,\Lambda_{i}\,P_{\,V_{i}}\,(\,f\,) \,\right\|_{X}^{\,p}$\, and \,$\sum\limits_{\,j \,\in\, J}\, w_{j}^{\,p}\, \left\|\,\Gamma_{j}\,P_{\,W_{j}}\,(\,g\,) \,\right\|_{Y}^{\,p}$\, are non-zero.\,Then
\begin{align*}
&\dfrac{A\,\left\|\,g \,\right\|_{Y}}{\left(\sum\limits_{\,j \,\in\, J}\, w_{j}^{\,p}\, \left\|\,\Gamma_{j}\,P_{\,W_{j}}\,(\,g\,) \,\right\|_{Y}^{\,p}\right)^{1 \,/\, p}}\,\left\|\,f \,\right\|_{X} \,\leq\, \left(\,\sum\limits_{\,i \,\in\, I} v_{i}^{\,p}\left\|\,\Lambda_{i}\,P_{\,V_{i}}\,(\,f\,) \,\right\|_{X}^{\,p}\,\right)^{1 \,/\, p}\\
&\hspace{3cm} \,\leq\, \dfrac{B\,\left\|\,g \,\right\|_{Y}}{\left(\sum\limits_{\,j \,\in\, J}\, w_{j}^{\,p}\, \left\|\,\Gamma_{j}\,P_{\,W_{j}}\,(\,g\,) \,\right\|_{Y}^{\,p}\right)^{1 \,/\, p}}\,\left\|\,f \,\right\|_{X}\\
&\Rightarrow\, A_{1} \,\left \|\,f \,\right \|_{X} \,\leq\, \left(\,\sum\limits_{\,j \,\in\, J}\,v_{j}^{\,p}\,\left\|\,\Lambda_{j}\,P_{\,V_{j}}\,(\,f\,) \,\right\|_{X}^{\,p}\,\right)^{1 \,/\, p} \,\leq\, B_{1} \,\left\|\,f\, \right\|_{X}\; \;\forall\; f \,\in\, X,
\end{align*} 
where 
\[A_{1} \,=\, \min\limits_{g \,\in\, Y}\left\{\dfrac{A\,\left\|\,g \,\right\|_{Y}}{\left(\sum\limits_{\,j \,\in\, J}\, w_{j}^{\,p}\, \left\|\,\Gamma_{j}\,P_{\,W_{j}}\,(\,g\,) \,\right\|_{Y}^{\,p}\right)^{1 \,/\, p}}\,\right\}\] and
\[B_{1} \,=\, \max\limits_{g \,\in\, Y}\left\{\dfrac{B\,\left\|\,g \,\right\|_{Y}}{\left(\sum\limits_{\,j \,\in\, J}\, w_{j}^{\,p}\, \left\|\,\Gamma_{j}\,P_{\,W_{j}}\,(\,g\,) \,\right\|_{Y}^{\,p}\right)^{1 \,/\, p}}\,\right\}.\]\;This shows that \,$\Lambda$\, is a \,$g$-$p$-fusion frame for \,$X$\, with respect to \,$\left\{\,X_{i}\,\right\}_{\, i \,\in\, I}$.\;Similarly, it can be shown that \,$\Gamma$\, is \,$g$-$p$-fusion frame for \,$Y$\, with respect to \,$\left\{\,Y_{j}\,\right\}_{\, j \,\in\, J}$.\\

Conversely, suppose that \,$\Lambda$\, and \,$\Gamma$\, are \,$g$-$p$-fusion frames for \,$X$\, and \,$Y$.\;Then there exist positive constants \,$A,\, B$\, and \,$C,\, D$\, such that
\begin{equation}\label{eqp1.4}
A \,\left\|\,f \,\right\|_{X} \,\leq\, \left(\,\sum\limits_{\,i \,\in\, I}\, v_{i}^{\,p} \,\left\|\,\Lambda_{i}\,P_{\,V_{i}}\,(\,f\,) \,\right\|_{X}^{\,p}\,\right)^{1 \,/\, p}  \,\leq\, B\, \left\|\, f \, \right\|_{X}\; \;\forall\; f \,\in\, X
\end{equation}
\begin{equation}\label{eqp1.5}
C\,\left\|\,g \,\right\|_{Y} \,\leq\, \left(\,\sum\limits_{\,j \,\in\, J}\, w_{j}^{\,p}\, \left\|\,\Gamma_{j}\,P_{\,W_{j}}\,(\,g\,) \,\right\|_{Y}^{\,p}\,\right)^{1 \,/\, p} \,\leq\, D\,\left\|\, g \, \right\|_{Y}\; \;\forall\; g \,\in\, Y.
\end{equation}
Multiplying (\ref{eqp1.4}) and (\ref{eqp1.5}), and using (\ref{eqp1.3}), we get
\[ AC\left\|\,f \,\otimes\, g\,\right\| \leq \left(\sum\limits_{i,\, j}v^{\,p}_{\,i}\,w^{\,p}_{\,j}\left\|\,\Lambda_{i}\,P_{\,V_{i}}\,(\,f\,) \otimes\, \Gamma_{j}\,P_{\,W_{j}}\,(\,g\,) \,\right\|^{\,p}\right)^{1 \,/\, p} \leq B\,D \left\|\,f \,\otimes\, g\,\right\|.\]
Therefore, for each \,$f \,\otimes\, g \,\in\, H \,\otimes\, K$, we get
\[ AC\left\|\,f \otimes g\,\right\| \leq \left(\sum\limits_{i,\, j}v^{\,p}_{\,i}\,w^{\,p}_{\,j}\left\|\,\left(\,\Lambda_{i} \otimes \Gamma_{j}\,\right)\,P_{\,V_{i} \,\otimes\, W_{j}}\,(\,f \otimes g\,)\,\right\|^{\,p}\right)^{1 \,/\, p} \leq BD\left\|\,f \otimes g\,\right\|.\]
Hence, \,$\Lambda \,\otimes\, \Gamma$\, is a \,$g$-$p$-fusion frame for \,$X \,\otimes\, Y$\, with respect to \,$\left\{\,X_{i} \,\otimes\, Y_{j}\,\right\}_{\,i,\,j}$\, with bounds \,$A\,C$\, and \,$B\,D$.\,This completes the proof. 
\end{proof}

\end{document}